\documentclass[12pt]{article}
\usepackage{latexsym}
\usepackage{amsfonts}
\makeatletter

\def\@footnotetext#1{\insert\footins{%
\footnotesize
    \interlinepenalty\interfootnotelinepenalty
    \splittopskip\footnotesep
    \splitmaxdepth \dp\strutbox \floatingpenalty \@MM
    \hsize\columnwidth \@parboxrestore
   \edef\@currentlabel{\csname p@footnote\endcsname\@thefnmark}\@makefntext
    {\rule{\z@}{\footnotesep}\ignorespaces
      #1\strut}}}
\def\abstract{\small\quotation{\hskip-\parindent\sc Abstract.}}

\def\classification{\@ifnextchar [{\@xfootnotenext}%
   {\begingroup\let\protect\noexpand
      \xdef\@thefnmark{}\endgroup
    \@footnotetext}}

\title {}

\begin{document}
\classification {{\it 2000 Mathematics
Subject Classification:} Primary 20F36,  secondary 20E36, 57M05.}

\begin{center}
{\bf \large  Representing braids by automorphisms} 
\bigskip

{\bf  
  Vladimir Shpilrain}

\end{center} 
\medskip

\begin{abstract}
\noindent 
 Based on a normal form for  braid group elements 
  suggested by Dehornoy, we 
prove several representations of braid groups by automorphisms of 
a free group to be faithful. This includes a simple proof of the 
standard Artin's representation being faithful. 
\end{abstract}

\date{}

\bigskip

\section{Introduction}

 Braid groups need no introduction; we just refer to the monograph 
\cite{B} for the background. Some notation has to be  reminded  
though. We denote the braid group on $n$ strands by $B_n$; this group
has a standard presentation $\langle \sigma_1, ..., \sigma_{n-1}| ~
\sigma\sb i \sigma\sb j = \sigma\sb j \sigma\sb i  ~{\rm if}
 ~\vert i-j\vert  > 1;  
~\sigma\sb i \sigma\sb {i+1} \sigma\sb i = \sigma\sb {i+1} \sigma\sb i \sigma\sb {i+1} 
 ~{\rm for}\ 1 \le i \le n-2 \rangle$.
We shall call elements of $B_n$ {\it braids}, as opposed to {\it braid words} that are
elements of the ambient free group on $\sigma_1, ..., \sigma_{n-1}$. We say that two 
braid words are {\it equivalent} if they represent the same braid. 

 There is a well-known representation (due to Artin) of the group $B_n$ in the
group  $Aut(F_n)$ of automorphisms of the free group $F_n$ (see e.g. 
\cite[p.25]{B}). Let $F_n$ be 
generated by $x_1, ..., x_n$. Then the automorphism $\hat \sigma_i$ 
corresponding to the braid generator $\sigma_i$, takes  $x_i$ to 
$x_i x_{i+1} x_i^{-1}$,  ~$x_{i+1}$ to $x_i$, and fixes all other free
generators.

 More recently, Wada \cite{W} has discovered several other representations 
of  the group $B_n$  by automorphisms of $F_n$. Some of them are obviously 
non-faithful; two of the remaining 4 are conjugate, which leaves us with the 
  following 3 interesting  representations: 
\medskip

\noindent {\bf (1)} This is actually an infinite series of representations
generalizing the standard Artin's representation. For an arbitrary non-zero integer 
$k$, the automorphism $\hat \sigma^{(k)}_i$ corresponding to the braid 
generator $\sigma_i$, takes  $x_i$ to 
$x_i^k x_{i+1} x_i^{-k}$,  ~$x_{i+1}$ to $x_i$, and fixes all other free
generators. 

\medskip

\noindent {\bf (2)} Here the automorphism $\hat \sigma_i$ 
corresponding to the braid generator  $\sigma_i$, takes  $x_i$ to 
$x_i x_{i+1}^{-1} x_i$,  ~$x_{i+1}$ to $x_i$, and fixes all other free
generators. 

\medskip

\noindent {\bf (3)} Here   $\hat \sigma_i$ takes  $x_i$ to 
$x_i^2 x_{i+1}$,  ~$x_{i+1}$ to $x_{i+1}^{-1} x_i^{-1} x_{i+1}$, and fixes   
all other free generators. 
\medskip

 In this paper, we prove the following  
\medskip 

\noindent {\bf Theorem A.} Each of the representations (1), (2), (3) 
above is faithful. 
\medskip 

 Our (very easy) proof is based on part (a) of  the following theorem 
 of Dehornoy \cite{P1}, \cite{P2}. At the same time,  our  proof 
establishes  part (b) of Dehornoy's theorem. In fact, each of the 
 faithfulness results of Theorem A 
 gives  a simple proof of part (b) of Theorem B. In particular, the
argument involving Wada's representation (3) seems to give
the   easiest proof of (b) known so far.

We call a braid word $u$ ~$\sigma_1$-nonnegative if there are no occurences 
of $\sigma_1^{-1}$ in $u$,  and $\sigma_1$-negative if there are no 
occurences 
of $\sigma_1$ with a positive exponent  in $u$. Then:
\medskip 

\noindent {\bf Theorem B.} (Dehornoy \cite{P1}, \cite{P2}) 

\noindent {\bf (a)} Every braid word is equivalent to either a 
$\sigma_1$-nonnegative  or a $\sigma_1$-negative braid word. 

\noindent {\bf (b)} If $u$ is a $\sigma_1$-nonnegative braid word 
with at least one occurence of $\sigma_1$, then $u$ is not equivalent to the 
empty word. 
\medskip 

 It is hoped that Dehornoy's normal form for  braid group elements
can be useful in proving   other representations of braid groups to be
faithful. We note at this point that part of our Theorem A follows 
from \cite[Theorem 7]{BH}, because it is proved there that, if $N_m$
 is the normal closure of the $n$ elements $x_1^m, \dots, x_n^m$ in $F_n$,
 where $m \ge 2$, then the induced action of the 
standard Artin's representation
on $F_n/N_m$  is faithful. Therefore, if a representation 
$\phi : B_n \to Aut(F_n)$   induces the same action on $F_n/N_m$ for some $m$
 as Artin's representation does (which  requires, in particular,  
$N_m$ being  invariant under $\phi(B_n)$), 
then  this $\phi$ must be faithful, too. Using this result, one
can establish faithfulness of, say, representations (1) above for $k \ne 2$.
 However, 
the combination of the  two conditions   ($N_m$ being  invariant under $\phi(B_n)$ 
and $\phi$  inducing the same action as Artin's representation on $F_n/N_m$) 
appears  to be
rather  restrictive, and is unlikely to be satisfied by most representations. 
For example, the representation (3) above satisfies the former condition 
for $m = 2$, but does not satisfy the latter. 
 Our method based on Dehornoy's normal form therefore appears to be more 
flexible. The referee has pointed out   that 
Larue  \cite{Larue} has used  a method similar  to ours    to show that 
the standard Artin's representation is faithful.

\smallskip 

 In the concluding Section 3, we show that  different Wada's representations 
have  different   images in $Aut(F_n)$,  with   one possible exception. A
probably difficult question is whether  or not different Wada's 
representations are conjugate. For example, 
 take two representations   $\varphi$ and $\psi$ 
 of type (1), where $\varphi: \sigma_i \to \hat \sigma^{(k)}_i; 
~\psi: \sigma_i \to \hat \sigma^{(-k)}_i$ ~for some non-zero integer $k$.  
Then  the images $\varphi(B_n)$ and $\psi(B_n)$  are conjugate
by the automorphism that takes every free generator $x_i$ to 
its inverse. This might be the only instance of different representations 
of the types (1)--(3) being conjugate, but I was not able to prove that. \\

\section{Proof of Theorem A  }
\bigskip

  Theorem A will  be proved if we establish  the following 
\medskip 

\noindent {\bf Lemma.} Let $\sigma \to \hat \sigma$ be  any of Wada's 
representations. Suppose $\sigma$ is a $\sigma_1$-positive braid word of the 
form  $\sigma_1 \sigma'$. 
Then  $\hat \sigma(x_1)$  has at least 2 occurences of $x_1^{\pm 1}$.  
\medskip 

\noindent {\bf Proof. ~(1)} We start with Wada's representation of type (1).
Since we assume that $\sigma$ is of the 
form  $\sigma_1 \sigma'$, we have the automorphism $\hat \sigma^{(k)}_1$ 
applied first, hence $\hat \sigma^{(k)}_1(x_1)=x_1^k x_2 x_1^{-k}$ 
already has at least  2  occurences of $x_1^{\pm 1}$.
        
 Then, any $\hat \sigma^{(k)}_i$ with $i \ge 2$ does not change existing 
occurences of $x_1$ and does not introduce any new ones.
 Thus, we have to only  concern ourselves with how $\hat \sigma^{(k)}_1$ 
acts on an element of the free group of the form $w=x_1^k u x_1^{-k}$, where 
 $u$ neither starts nor ends with $x_1^{\pm 1}$. 
We are going to show that $\hat \sigma^{(k)}_1(w)$ has the same form 
(with different $u$, perhaps), i.e., that $x_1^k$ on the left and $x_1^{-k}$
on the right cannot cancel after $\hat \sigma^{(k)}_1$ is applied. Because of 
the symmetry, we are going to consider $x_1^k$ on the left only. 
Consider 2 cases: 
\medskip 
        
\noindent {\bf (a)} $u=x_m^s u'$, where $m \ge 3, ~s \ne 0$,  and $u'$ does 
not  start with $x_m^{\pm 1}$. 
Then  $\hat \sigma^{(k)}_1(w)=x_1^k x_2^k x_1^{-k} x_m^s  
\hat \sigma^{(k)}_1(u') x_1^k x_2^{-k} x_1^{-k}$,  ~and  $x_1^k$ 
on the left does not cancel.  Indeed, for a cancellation process to start, 
there must be a cancellation between $x_m^s$  and  
$\hat \sigma^{(k)}_1(u')$, i.e., $\hat \sigma^{(k)}_1(u')$ should start 
with $x_m^{\pm 1}$. Since $u'$ itself does {\it not} start 
with $x_m^{\pm 1}$ and $\hat \sigma^{(k)}_1$ does   not affect occurences 
of $x_m$, that could only mean that some initial fragment of $u'$ 
 became the empty word  after $\hat \sigma^{(k)}_1$ was applied. But this
is impossible because $\hat \sigma^{(k)}_1$ is an automorphism. 
\medskip 
        
\noindent {\bf (b)} $u=x_2^s u'$, ~where $s \ne 0$, and  $u'$   does not 
start  with $x_2^{\pm 1}$.  Then $\hat \sigma^{(k)}_1(w)=x_1^k x_2^k x_1^{-k} x_1^s 
\hat \sigma^{(k)}_1(u')$, 
where there is no  cancellation between $x_1^s$ and  $\hat 
\sigma^{(k)}_1(u')$ 
because $\hat \sigma^{(k)}_1(u')$ cannot start with $x_1^{\pm 1}$. 
To have the cancellation process get to $x_2^k$, we must have $k=s$. 
  Then, to cancel all of $x_2^k$, we must have $\hat \sigma^{(k)}_1(u')$ 
start with $x_2^{-k}$, which is impossible. Indeed, if $u'$ starts   
with $x_m^{\pm 1}$, ~$m \ge 3$, then, in order to have $\hat 
\sigma^{(k)}_1(u')$ 
start with $x_2^{\pm 1}$, we must have some non-empty fragment of $u'$ 
between $x_m^{\pm 1}$ and $x_m^{\mp 1}$  become the empty word  after $\hat
\sigma^{(k)}_1$ is applied, which is impossible because $\hat \sigma^{(k)}_1$ 
is an automorphism. If $u'$ starts  with $x_1^{\pm 1}$, then the obvious 
inductive argument implies that  $\hat \sigma^{(k)}_1(u')$ should start with
$x_1^{\pm 1}$ as well. (Note that the length of $u'$ is smaller than that 
of $w$). 
        
 Thus, in either case,  $x_2^k$ cannot cancel, and therefore, $\hat
\sigma^{(k)}_1(w)$ has the  same form as $w$.  $\Box$ 
        
\medskip 
        
\noindent {\bf (2)} For Wada's representation of type (2), the proof 
goes along exactly the same lines. 
        
\medskip 
        
\noindent {\bf (3)} Finally, consider Wada's representation of type (3).
Again,  we have the automorphism $\hat \sigma_1$ 
applied first, hence $\hat \sigma_1(x_1)=x_1^2 x_2$ 
already has 2  occurences of $x_1$. 
        
 Also,  any $\hat \sigma_i$ with $i \ge 2$ does not change existing 
occurences of $x_1$ and does not introduce any new ones, so 
 we have to only  concern ourselves with how $\hat \sigma_1$ 
acts on an element of the free group of the form $w=x_1^2 u$, where 
 $u$  does not start with $x_1^{\pm 1}$. We are going to show that 
$\hat \sigma_1(w)$ has the same form. Again, there are 2 cases: 
\medskip 
        
\noindent {\bf (a)} $u=x_m^{\pm 1} u'$, where $m \ge 3$. Then 
$\hat \sigma_1(w)=x_1^2x_2 x_m^{\pm 1} \hat \sigma_1(u')$, hence 
 $x_1^2$ on the left does not cancel. (If $u'$ does not start with $x_m^{\mp 1}$,
 then neither does $\hat \sigma_1(u')$. Therefore, for   cancellation   
 between $u'$ and $\hat \sigma_1(u')$ 
 to  occur, some fragment of $u'$ must be mapped to the empty word, which is 
 impossible  since $\hat \sigma_1$ 
is an automorphism).
  
\medskip 
        
\noindent {\bf (b)} $u=x_2^s u'$, ~where $s \ne 0$, and  $u'$   does not 
start  with $x_2^{\pm 1}$. Then 
$\hat \sigma_1(w)=x_1^2x_2 x_1^2x_2 x_2^{-1}x_1^{-s} x_2 \hat \sigma_1(u')$, 
 and there is no cancellation between $x_2$ and $\hat \sigma_1(u')$ 
since $\hat \sigma_1(u')$ cannot start with $x_2^{-1}$ unless $u'$ 
starts  with $x_2^{\pm 1}$.
        
 Therefore, $\hat \sigma_1(w)=x_1^2x_2 x_1^{2-s} x_2 \hat \sigma_1(u')$, 
 and no matter what $s$ is, $\hat \sigma_1(w)$ has the form $x_1^2 u$, 
where  $u$  does not start with $x_1^{\pm 1}$.  $\Box$ \\

\section{Images of Wada's representations }
\bigskip

 Here we prove the following
\medskip 

\noindent {\bf Proposition.} Let $\varphi$ and $\psi$ be two 
different Wada's representations. Then the groups $\varphi(B_n)$ and
$\psi(B_n)$  are different subgroups of   $Aut(F_n)$ unless, perhaps,
  $\varphi$
and $\psi$  are both of type (1) and $\varphi: \sigma_i \to \hat 
\sigma^{(k)}_i; 
~\psi: \sigma_i \to \hat \sigma^{(-k)}_i$ ~for some non-zero integer $k$. 

\medskip 

\noindent {\bf Proof.} We have to consider several cases. 
\smallskip

\noindent {\bf (1)} Both $\varphi$ and $\psi$ are of type (1), so that
$\varphi: \sigma_i \to \hat \sigma^{(k)}_i; 
~\psi:  \sigma_i \to \hat \sigma^{(s)}_i$ ~for some non-zero integers $k, s$.
 We assume that $k \ne \pm s$, so let $|k| > |s|$.

 Consider the Magnus representation of the groups $\varphi(B_n)$ and
$\psi(B_n)$ ~(for $k=1$, it is also known as the Burau representation --
see \cite[p.102]{B}). Under this representation, the automorphism 
$\hat \sigma^{(k)}_i$ is mapped onto the $n \times n$ matrix which   
differs from the identity matrix only by a $2 \times 2$ block with 
the top left corner in the $(i, i)$th place. This block is 
$\left(\begin{array}{cc} 1-t^k & t^k \\ 1 & 0\end{array}\right)$. 
Thus, the determinant of this matrix is $-t^k$,  and therefore, 
the determinant of the matrix corresponding to an arbitrary braid under 
the composition of $\varphi$ and the Magnus representation, is equal 
to $\pm t^{km}$ for some integer $m$. 

 Now if we take, say, $\psi(\sigma_1)$ and then apply the Magnus 
representation,
we shall get a matrix with the determinant $-t^s$. A matrix like that 
cannot be a product of matrices with   determinants of the form $\pm t^{km}$ 
since $|k| > |s|$. This completes the proof in case (1).
\smallskip

\noindent {\bf (2)} $\varphi$ is of type (1), and $\psi$ is of type (2). 
Again, we apply the Magnus representation to both groups $\varphi(B_n)$ and
$\psi(B_n)$. Note that this is possible since the {\it mapping group} 
of $\psi(B_n)$ is the same as that of $\varphi(B_n)$. (The mapping group 
of a single automorphism $\alpha : x_i \to y_i, ~1 \le i \le n,$ is the group
with the presentation $\langle x_1, ..., x_n| ~x_i=y_i, ~1 \le i \le n 
\rangle$.
The mapping group of a group $G$ of automorphisms has the set of relations 
which
is the union of sets of relations for the  mapping group of 
each individual automorphism in $G$). 

 If we apply the Magnus representation to $\psi(\sigma_i)$, we shall get a 
matrix which  
differs from the identity matrix only by the following $2 \times 2$ block
with  the top left corner in the $(i, i)$th place:  
$\left(\begin{array}{cc} 2 & -1 \\ 1 & 0\end{array}\right)$. The determinant 
of this matrix is 1. Therefore, any matrix in the image of the  Magnus
representation of $\psi(B_n)$ has determinant 1. This completes the proof in 
case
(2).

\smallskip

\noindent {\bf (3)} $\varphi$ is of type (1) or (2), and $\psi$ is of type 
(3).
 In that case,  $\varphi(B_n)$ and
$\psi(B_n)$  are different subgroups of   $Aut(F_n)$ because they have 
different presentations of their mapping groups. More accurately, 
the mapping group of $\varphi(B_n)$ has the presentation 
$\langle x_1, ..., x_n| ~x_1=x_2=...=x_n \rangle$, whereas $\psi(B_n)$ has the
presentation $\langle x_1, ..., x_n| ~x_1=x_2^{-1}=...=x_n^{(-1)^{n+1}}
\rangle$. The groups themselves are isomorphic, yet the presentations 
are different. Now we argue as follows. 

  Consider, say, $\varphi(\sigma_1)$. The mapping group of this automorphism
has the presentation $\langle x_1, ..., x_n| ~x_1=x_2 \rangle$. 
Suppose $\varphi(\sigma_1) \in \psi(B_n)$. That means, in particular, that 
by adding some elements to  ${\{}x_1x_2^{-1}{\}}$, we
can get a set of elements of a free group whose normal closure is the same 
as that of  ${\{}x_1x_2, ...,  x_1x_n{\}}$. But this is  impossible
since the normal closure of  the union of these two sets contains, say, the 
element $x_1^2$, which none of the two normal closures alone does.  $\Box$ \\

\medskip
\noindent 
 Department of Mathematics\\ 
The City  College  of New York\\ 
New York, NY 10031 \\
 
\smallskip 

\noindent {\it e-mail address\/}: shpil@groups.sci.ccny.cuny.edu 

\smallskip

\noindent {\it http://zebra.sci.ccny.cuny.edu/web/shpil} 

\end{document}